\documentclass[MSNbibl,number,citesort,seceqn,dvips]{arxbj}

\aid{0}
\volume{19}
\issue{4}
\pubyear{2013}
\firstpage{1122}
\lastpage{1149}
\doi{10.3150/12-BEJSP04} 

\makeatletter
\newproclaim{GP}{General Principle}
\newproclaim{OP}{Open Problem}
\newtheorem{Proposition}{Proposition}[section]
\newtheorem{Theorem}{Theorem}[section]
\newtheorem{Conjecture}{Conjecture}
\newremark{Completegraph}{Complete graph}
\newremark{dimensionaldiscretetorus}{$d$-dimensional discrete torus $\Ints^d_m$}
\newremark{Smallworlds}{Small worlds}
\newremark{Randomgraphs}{Random graphs with prescribed degree distribution}
\newremark{Longrange}{Long range geometric interaction}
\newremark{Proximitygraphs}{Proximity graphs}
\newremark{TheHamming}{The Hamming cube $\{0,1\}^d$}
\newremark{TheTokenprocess}{The Token process}
\newremark{Anuance}{A nuance}
\newremark{ThePandemicprocess}{The Pandemic process}
\newremark{TheAveragingprocess}{The Averaging process}
\newremark{Votermodel}{Voter model}
\newremark{Coalescingchainsmodel}{Coalescing chains model}
\newremark{Localtransienceprinciple}{Local transience principle}
\newremark{TheCompulsiveGamblerprocess}{The Compulsive Gambler process}
\newremark{TheInterchangeprocess}{The Interchange process}
\newremark{TheDeferenceprocess}{The Deference process}
\newremark{TheFashionistamodel}{The Fashionista model}

\newcommand{\bX}{\mathbf{X}}
\newcommand{\eps}{\varepsilon}
\newcommand{\Agents}{\mathbf{Agents}}
\newcommand{\States}{\mathbf{States}}
\newcommand{\Reals}{{\mathbb{R}}}
\newcommand{\Ints}{{\mathbb{Z}}}
\newcommand{\CC}{{\mathcal C}}
\newcommand{\EE}{{\mathcal E}}
\newcommand{\FF}{{\mathcal F}}
\newcommand{\NN}{{\mathcal N}}
\renewcommand{\SS}{{\mathcal S}}
\newcommand{\VV}{{\mathcal V}}
\newcommand{\taumeet}{\tau_{\mathrm{meet}}}
\newcommand{\taumix}{\tau_{\mathrm{mix}}}
\newcommand{\cd}{\,\mathop{\to}\limits^{d}\,}
\newcommand{\ed}{\,\mathop{=}\limits^{d}\,}
\newcommand{\Ex}{{\mathbb E}}
\renewcommand{\Pr}{{\mathbb{P}}}
\newcommand{\Ent}{\operatorname{Ent}}
\newcommand{\Var}{\operatorname{var}}
\makeatother

\begin{document}
\begin{frontmatter}

\title{Interacting particle systems as stochastic social dynamics}
\runtitle{IPS as social dynamics}

\begin{aug}
\author{\fnms{David} \snm{Aldous}\corref{}\ead[label=e1]{aldous@stat.berkeley.edu}\ead[label=u1,url]{www.stat.berkeley.edu/\textasciitilde aldous/}}
\runauthor{D. Aldous} 
\address{Statistics Dept., U.C. Berkeley, CA 94720, USA.
\printead{e1};\\
\printead{u1}}
\end{aug}


%
\begin{abstract}
The style of mathematical models known to probabilists as
Interacting Particle Systems and exemplified by the
Voter, Exclusion and Contact processes
have found use in many academic disciplines.
In many such disciplines the underlying
conceptual picture is of a social network, where individuals
meet pairwise and update their ``state'' (opinion, activity etc) in a way
depending on the two previous states.
This picture motivates a precise general setup we call
\emph{Finite Markov Information Exchange (FMIE) processes}.
We briefly describe a few less familiar models
(Averaging, Compulsive Gambler, Deference, Fashionista)
suggested by the social network picture,
as well as a few familiar ones.
\end{abstract}

%

\end{frontmatter}

\section{Introduction}

What is the most broad-ranging currently active field of applied
probability? Measuring breadth by the number of different academic
disciplines where it appears, it seems hard to beat a field appearing
in Physics, Computer Science and Electrical Engineering, Economics and
Finance, Psychology and Sociology, Population Genetics and Epidemiology
and Ecology. Curiously, the field has no good name, though readers of
\emph{Bernoulli} will most likely know it, under the
continuing influence of Liggett's 1985 book~\cite{Liggett},
as \emph{Interacting Particle Systems}.
It turns out that mathematically similar toy models of different
real-world entities have been repeatedly re-invented in the different
disciplines mentioned above, and literally thousands of papers have
been written
since 2000 on discipline-dependent
variant models intended as more realistic.

Rather than attempting a brief overview of the whole field as it stands today
(on which I am certainly not an expert), let me take an opposite
approach, imagining starting
with a blank slate. I~will set out in the next section a
particular conceptual setting -- what sorts of things are we trying to
model? -- then a particular technical setup, and then indicate what
sort of
mathematical results one might seek.
Given this manifesto, the rest of the article describes some basic models,
gives pointers to the existing literature and presents a few open
problems suggested by our setting.
The article is based on a short course given in summer 2012 at Warwick
and Cornell, and more extensive informal material is available at \cite
{aldous-FMIE}. I must emphasize that the article is intended as an
introduction to the topics it mentions, not as an
authoritative survey,
which would require ten times the length!
In particular none of the mathematics outlined here is both new and
non-straightforward, except perhaps some open problems.

\section{Conceptual and technical setup}

\subsection{Information flow through social networks via pairwise interaction}

The particular setting we envisage is of the human social network.
People meet or communicate, and explicitly or implicitly receive
information from the other, and this interaction may change the
``state''
(opinion or possessions or \ldots) of the individuals.
Now that sentence skips over many issues.
``Information'' is a complicated concept~\cite{Floridi},
and formalizing it as one of a set of prespecified states hardly
captures the nuances of human experience. And obviously we receive
information in many ways other than
one-on-one interaction with other individuals.
But the fact that humans do a lot of different things via interactions
emphasizes the point that we are interested in a variety of different
rules, a point
that would be less visible if we took
some real or hypothetical physical
system as a prototype setting.

Social networks (whose mathematical modeling and analysis long predates
the Facebook era) are usually modeled as unweighted
graphs -- vertices are individuals, edges indicate some specified relationship.
In most contexts it would be more realistic to use \emph{weighted}
graphs, the
weight indicating strength of relationship.
This is seldom done because of increased model complexity and the fact
that the meaning of ``strength of relationship'' is context-dependent.
But we will give it the specific meaning of ``frequency of meeting'',
formalized via the rates $\nu_{ij}$ below.

The story above suggests the following technical framework, in
which randomness arises primarily via assumed random times of meetings.
Because the framework can be used outside the social network setting,
I will use the impersonal word ``agent'' in what follows.

\subsection{The technical framework}
\label{sec:tech}

Consider $n$ \emph{agents}
and a nonnegative array
$(\nu_{ij})$, indexed by unordered pairs $\{i,j\}$,
which is irreducible
(i.e., the graph of edges corresponding to strictly positive entries is
connected).
Assume
\begin{itemize}
\item Each unordered pair $i,j$ of agents with $\nu_{ij} > 0$ meets at
the times of a rate-$\nu_{ij}$ Poisson process, independent for
different pairs.
\end{itemize}
Call this collection of Poisson processes the \emph{meeting process}.
Where needed we write $\Agents$ for the set of agents.

Now introduce a set $\States$ of states, and consider a process
$\bX(t)$ in which each agent $i$ is in some state
$X_i(t)$ at time $t$, and the state can only change when agent $i$
meets some other agent $j$, at which time the two agents' states are
updated according to some rule (typically deterministic, but might be
random) which
depends only on the pre-meeting states $X_i(t-), X_j(t-)$.
So in the deterministic case, there is a given arbitrary function
$F: \Agents\times\Agents\to\Agents$ and when agents in states
$s_1$ and $s_2$ meet their states are updated as
%
%
\begin{equation}
\label{def-F} (s_1,s_2) \to \bigl(F(s_1,s_2),
F(s_2,s_1) \bigr).
\end{equation}
We call such a process an \emph{FMIE (Finite Markov Information Exchange)
process} and emphasize that a specific such process has two levels.
The meeting process is determined by the meeting rate matrix $(\nu_{ij})$
(which we loosely call the ``geometry'', imagining that we tend to meet
more frequently with individuals who are ``close'' in some sense).
The information update process is determined by some rule.
Then the different named FMIE models we shall describe (such as the
Voter model) correspond to different rules for the informational
superstructure; we regard such models as defined over
arbitrary geometries.

\subsection{Three disciplines using FMIE models}

Of course much of the beauty and utility of mathematics is in providing
abstractions that encompass apparently different things, and as stated
at the start FMIE models have been used in many disciplines. Let me
briefly mention three that are relevant to our models later.

Conceptually closest is the mathematical modeling of epidemics,
whose long history is indicated by the 1957 monograph~\cite{MR0095085}.
A typical classical model like the SIR
(Susceptible--Infected--Recovered) model corresponds to a 3-state
FMIE if the relevant distributions are Exponential.\footnote
{Capitalizing named distributions is a useful convention that should be
more widely adopted.}
Here of course the vertices are people or animals and the edges
relate to physical proximity.
Classical models considered homogeneously mixing populations
(i.e., the complete graph geometry) but there is extensive post-2000
literature on more general geometries.

A quite different context, from theoretical computer science, concerns
distributed randomized algorithms. Here agents are processors,
communicating with nearby (direct connections in a network) processors,
and performing computations based on data transferred. See~\cite{Shah}
for the particular topic of \emph{gossip algorithms}.

What are called \emph{Interacting Particle Systems (IPS)} arose historically
as statistical physics models of phenomena such as magnetism, and mathematicians
became interested in particular in rigorous study of phase transitions,
in the setting of
infinite geometries, specifically the infinite $d$-dimensional lattice.
This literature long ago branched out to study many ``toy models'' not
so directly related to any specific physical phenomenon.
Our development is in the same ``toy models'' spirit, and the models
naturally suggested
by the social networks context partly overlap and partly differ from
those naturally suggested by statistical physics.
In particular, the abstraction of what we are studying as ``information
flow through networks'' seems more helpful than any physical visualization.

A broader but less focused view of the whole ``statistical physics of
social dynamics'' topic can be found in
\cite{castellano}.

\subsection{Style of mathematical results}

The concept of an \emph{IPS} is generally discussed via example rather
than explicit definition. For instance Liggett~\cite{Liggett} takes as examples
the Voter model, Contact process, Exclusion process, Glauber dynamics
for the Ising model, and the Potlatch process.
Our mathematical setup of \emph{FMIE process} in Section~\ref
{sec:tech} is
intended as a more precisely circumscribed class of processes to study.
For instance, only the first three processes on Liggett's list are FMIE
(for the contact process, add an artificial agent with rate-one
meetings with each other agent, to effect recovery).
What we call the meeting process is often called the
\emph{graphical representation} in IPS.

Our emphasis on a \emph{finite} number of agents gives a somewhat
different perspective than in classical IPS, in which the finite
setting was considered only occasionally, e.g.,~\cite{MR958203}.
Our viewpoint is by analogy with the
modern literature on mixing and hitting times for finite Markov chains
\cite{Levin,aldous-fill}; study how the behavior of a model depends on
the underlying geometry, where ``behavior'' means quantitative aspects
of finite-time behavior.
More precisely, we typically study $n \to\infty$ limits in a family
of finite-agent processes, rather than $t \to\infty$ limits in a
single infinite-agent process.
The Averaging process studied in Section~\ref{sec:Ave} provides a
prototype for the style of results we would like to establish for a
given FMIE model; see also the discussion of possible $t \to\infty$
limits in Section~\ref{sec:SCC}.

The FMIE setup allows an arbitrary update function $F$ at (\ref
{def-F}) and an arbitrary ``geometry'' of meeting rates $(\nu_{ij})$.
One can hardly hope to say anything of substance in this generality
(we will mention various ``general principles'' but these are rather
trite, once said).
At a pragmatic level, if seeking analytic results one chooses an update
function such that the FMIE process can be
analyzed on the complete graph (``mean-field'', roughly) geometry,
otherwise it will likely be difficult to handle other geometries.
And obversely, one chooses a geometry in which one can analyze the
two basic processes of Section~\ref{sec:2basic}, otherwise
it will likely be difficult to handle more complicated update rules.

From our viewpoint, much of the literature
since 2000 consists of inventing and investigating a large number of
increasingly elaborate update rules, while considering only a small
number of ``standard geometries'', indicated below.
Keep in mind that it is
easy to invent and simulate models, but hard to give rigorous proofs or
to relate convincingly to real-world data.\footnote{In talks I
criticize a particularly egregious instance~\cite{SNARC} of claims of
real-world significance based on models with absolutely no data or
evidence linking the model to the asserted real-world phenomena.}

\subsection{Some standard geometries}

Here we briefly describe four of the ``standard geometries''
on which models are often studied.
Recall $n$ is the number of agents.
Note that in a realistic social network there is no reason why
$\sum_j \nu_{ij}$, which is the rate at which agent $i$ meets other agents,
should be constant in $i$.
But many simple models
do have this property,
in which case we may rescale time to fit a
\emph{standardization convention}
%
%
\begin{equation}
\label{standardization} \sum_j
\nu_{ij} = 1\qquad \forall i
\end{equation}
as in the first two examples below.

\begin{Completegraph*}
Here $\nu_{ij} = 1/(n-1), j \ne i$.
\end{Completegraph*}

\begin{dimensionaldiscretetorus*}
Here $\nu_{ij} = 1/(2d)$ for adjacent agents $i, j$.
\end{dimensionaldiscretetorus*}

The next two examples involve networks modelled as \emph{random} graphs.
Formally this adds another level of randomness, though informally we
can just think of being given a typical realization.
In such (non-regular) graphs
one typically cannot satisfy (\ref{standardization})
and instead the default convention is to take rates $\nu_{ij} = 1$ for
each edge $(i,j)$.

\begin{Smallworlds*}
Here one starts with the $d$-dimensional discrete torus (or similar
non-random graph) and then adds random other edges in some specified
way, for instance
independently with probabilities proportional to
$(\mbox{Euclidean length})^{- \gamma}$.
\end{Smallworlds*}

\begin{Randomgraphs*}
There are several ways
(in particular the \emph{configuration model}~\cite{bollobas})
to define graphs corresponding to the intuitive notion of being
``completely random'' subject to the constraint that, as
$n \to\infty$, the distribution of vertex-degrees converges to a
prescribed distribution. Such models are tractable because their
$n \to\infty$ local weak limit is just a modified Galton-Watson tree.
\end{Randomgraphs*}

\subsection{Some less standard geometries}

There is a huge ``complex networks'' literature (see, e.g., the monographs
\cite{dorog,newman}) devoted to
inventing and studying network models,
and I will not write out a long list here.
Let me instead mention three other mathematically natural geometries which
have \emph{not} been studied as much as one might have expected
(and two more will be mentioned in Section~\ref{sec:infinite}).

\begin{Longrange*}
Start with the torus $\Ints^d_m$ and add rates
$\nu_{ij} = c_{d,m,\gamma} \Vert i - j\Vert_2^{-\gamma}$
for non-adjacent $i,j$.
This is a way to interpolate between the torus and the complete graph.
Note the distinction between this and the ``small worlds'' model, in
which the rate is
non-vanishing for a few random edges.
\end{Longrange*}

\begin{Proximitygraphs*}
Given a model whose behavior is understood on the two-dimensional lattice,
one could investigate the effect of ``disorder'' by taking instead
a proximity graph over a Poisson process of points in the plane \cite
{aldous-shun}.
\end{Proximitygraphs*}

\begin{TheHamming*}
This is a standard example in Markov chain theory.
\end{TheHamming*}

\section{The two background models}
\label{sec:2basic}

Here we discuss two FMIE processes that are ``basic'' in two senses.
They fit the definition of FMIE process but are much simpler than
(and therefore unrepresentative of) a typical FMIE process; yet many
other FMIE processes can be regarded as being built over one of these
base processes, with extra structure added.

For each FMIE, we give only a verbal description, which the reader can readily
formalize.

\subsection{The Token process and the associated Markov chain}

\begin{quote}
\begin{TheTokenprocess*}
There is one token. When the agent $i$ holding the token meets another
agent $j$, the token is passed to $j$.
\end{TheTokenprocess*}
\end{quote}
This can evidently be viewed as a 2-state FMIE process.
The natural aspect to study is $Z(t) = $ the agent holding the token at
time $t$.
This $Z(t)$ is the (continuous-time) Markov chain with transition rates
$(\nu_{ij})$.
We call this the \emph{associated} Markov chain\footnote{Throughout,
Markov chains are continuous-time.} -- associated with the underlying geometry.

\begin{Anuance*}
Let me digress to mention, by analogy, a nuance.
A Google Scholar search on
``exact phrase $=$ Galton-Watson;
year $=$ 1965--1969''
returns only papers studying
$Z_g = $ population size in generation $g$ of a Galton-Watson branching process.
Changing to ``year $=$ 2005--2009'', half the papers talk about the
Galton-Watson \textbf{tree}, which of course has more structure than the
population size process -- one can study cousins, for instance.
Analogously, the Token process has more structure than the
associated Markov chain -- one can study the time until
agent $i$ meets someone who has previously
held the token, for instance.
\end{Anuance*}

For finite-state Markov chains, we have available
both a classical general theory
and a long catalog of explicit calculations in specific geometries.
For our purposes, the main point is that several FMIE models
(classically the Voter model and the exclusion process, but also the
Averaging process) have close connections with the associated Markov
chain, and parts of the ``standard modern theory'' \cite
{Levin,aldous-fill} of finite chains can be used in their analysis.

For readers more familiar with discrete-time chains, let me emphasize that
(because $\nu_{ij} = \nu_{ji}$) the stationary distribution of the
associated chain is always uniform; this does not assume
the standardization (\ref{standardization}).

\subsection{The Pandemic process and FPP}
\label{sec:Pan}

\begin{quote}
\begin{ThePandemicprocess*}
Initially one agent is infected. Whenever an infected agent meets
another agent, the other agent becomes infected.
\end{ThePandemicprocess*}
\end{quote}
In classical language, this is the SI epidemic with $\operatorname{Exponential}(\nu_{ij})$ distributions for infection times.
But again there is a nuance.
Suppose that on each undirected edge $(i,j)$ of the geometry we put a single
$\operatorname{Exponential}(\nu_{ij})$ r.v. $\xi_{ij}$, independent over edges, and regard
$\xi_{ij}$ as the length of the edge.
Then there is an induced distance ($=$ minimum route length) between
any pair of agents.
This model of random distances is what I will call the
\emph{first passage percolation (FPP)} model.
Note this differs from classical lattice FPP, which assumes i.i.d.
variables with arbitrary distribution.

Here is the nuance.
In the Pandemic process started with agent $i$, write $T^{\mathrm
{pan}}_{ij}$ for
the time when agent $j$ is infected.
In the FPP process write $T^{\mathrm{fpp}}_{ij}$ for the distance
between $i$ and $j$.
From the memoryless property of the Exponential, we see that
$T^{\mathrm{pan}}_{ij} \ed T^{\mathrm{fpp}}_{ij}$ and indeed that
$(T^{\mathrm{pan}}_{ij}, j \in\Agents) \ed(T^{\mathrm{fpp}}_{ij},
j \in\Agents)$ for
each $i$.
But Pandemic and FPP are different processes in that
$T^{\mathrm{fpp}}_{ij} = T^{\mathrm{fpp}}_{ji}$ whereas we have only
$T^{\mathrm{pan}}_{ij} \ed T^{\mathrm{pan}}_{ji}$.
(One\vspace*{2pt} can fix this particular issue by re-interpreting, in FPP, the time
$\xi_{ij}$ as the elapsed time, after one end is infected, until the other
end is infected, but this misses the overall point that the distributions
of the whole arrays $(T_{ij})$ are different.)

The Pandemic model is ``basic'' in a very specific sense, in that it exhibits
the fastest possible spread of information in any FMIE model.
That is, in an arbitrary FMIE process $\bX(t)$ over a given geometry,
the value of $X_j(t)$ can only be influenced by the value of $X_i(0)$
if $T^{\mathrm{pan}}_{ij} \le t$.

Here of course we are exploiting the natural coupling of two arbitrary
FMIE processes over the same geometry, which simply uses the same
realization of the meeting process.\footnote{If both update rules are
random then we need also to specify some joint distribution.}
Albeit obvious, one can view this as an instance of the following.
\begin{GP}
It is often useful to consider the natural coupling of two
FMIE processes.
\end{GP}

Models of real-world epidemics are of course much more elaborate.
Regarding our simple Pandemic model, it is analogous to the Markov
chain in two ways:
\begin{itemize}
\item There is extensive literature doing calculations in specific geometries;
I show a little in Section~\ref{sec:Pandemic}.
\item Several other FMIEs can be ``built over'' the
Pandemic model -- see Section~\ref{sec:built_over_P}.\vadjust{\goodbreak}
\end{itemize}
But it differs it that there seems no existing ``general theory''. See
a conjecture in Section~\ref{sec:WLLN}.

\subsection{Remarks on time-asymptotics}
\label{sec:SCC}

In these remarks, we are envisaging the finite-state case.
\begin{GP}
Time-asymptotics for fixed $n$ are not the issue;
in most models it's easy to see what's going on.
\end{GP}

Informally there are four qualitatively different possible time-asymptotic
behaviors, and each of the models we treat fits readily into
one of these.
\begin{itemize}
\item Absorption in a random one (of a small number of) ``ordered''
configurations. [Pandemic, Averaging, Voter, Deference].
\item Absorption in a random one (of a large number of) ``disordered''
configurations. [Compulsive Gambler].
\item\ldots\ldots\ldots
\item Convergence to the unique stationary distribution (perhaps on a
subset of configurations). [Token, Interchange, Fashionista].
\end{itemize}

The elided third item is the most general possibility for a
finite-state Markov chain
(applied to the configurations of a FMIE process);
there are one or more \emph{strongly connected\footnote{Meaning there
is a directed path between any ordered pair of states.} components}
(SCCs) of the
directed transition graph, and the
limit distribution is a mixture over the stationary distributions on
each SCC.
While it's easy to invent artificial FMIE models with this behavior, I~don't know any ``natural'' such model.

The following question is less easy that it first appears.
The natural conjecture (involving subsets of $\Agents$ closed under $F$,
and functions whose sum is conserved by the interaction) is false.
\begin{OP}
Consider the FMIE process over the complete graph on $n$ agents (for
large $n$)
defined by a deterministic update function $F$ as at (\ref{def-F}),
What are the SCCs of the FMIE process?
\end{OP}

\section{A prototype FMIE model: The Averaging process}
\label{sec:Ave}

As a prototype FMIE model -- for which one can say something
quantitative about its behavior over a general geometry -- we
take the following.
Curiously, though it appears as an ingredient in more elaborate models,
the only place I have seen the basic process itself is in the ``gossip
algorithms'' literature~\cite{Shah} which derives a version of
the ``global bound'', Proposition~\ref{PGcon}.
A~careful treatment of the results below is written out in \cite
{aldous-lanoue} -- here I give a compressed version.

\begin{quote}
\begin{TheAveragingprocess*}
Initially each agent has some (real-valued) amount of money.
Whenever two agents meet, they divide their money equally.
\end{TheAveragingprocess*}
\end{quote}
It is rather obvious that the individual values $X_i(t)$ converge as $t
\to\infty$ to the average
of initial values, so what is there to say?

Write $\mathbf{1}_i$ for the initial configuration $(1_{(j = i)}, j
\in\Agents)$,
that is agent $i$ has unit money and other agents have none,
and write $p_{ij}(t)$ for the transition probabilities of the
associated Markov chain.
\begin{Proposition}
\label{L1}
For the Averaging process with initial configuration $\mathbf{1}_i$ we have
$\Ex X_j(t) = p_{ij}(t/2)$.
More generally, from any deterministic initial configuration $\mathbf
{x}(0)$, the expectations
$\mathbf{x}(t) := \Ex\bX(t)$ evolve exactly as the dynamical system
\[
\frac{d}{dt} \mathbf{x}(t) = {\frac{1}{2}} \mathbf{x}(t) \NN,
\]
where $\NN$ is the generator of the associated Markov chain.
\end{Proposition}
The key point behind this results is that we can rephrase the dynamics
of the Averaging process as
\[
\mbox{\emph{when two agents meet, each gives half their money to the other}.}
\]
In informal language, this implies that the motion of
a random penny -- which at a meeting of its owner agent is given to the
other agent with probability
$1/2$ -- is as the associated Markov chain at half speed, that is with
transition rates $\nu_{ij}/2$.

To continue, we need to
recall some background facts (e.g.,~\cite{aldous-fill}, Chap. 3) about
reversible chains, here specialized to our continuous-time setting of
the associated Markov chain,
that is with uniform stationary distribution ($\nu_{ij} = \nu_{ji}$).
A function $f \dvtx\Agents\rightarrow\Reals$ has
(with respect to the uniform distribution)
average $\overline{f}$, variance $\Var{f}$ and $L^2$ norm $\Vert
f\Vert_2$ defined by
\[
\overline{f} := n^{-1} \sum_i
f_i, \qquad \Vert f\Vert_2^2 :=
n^{-1} \sum_i f_i^2,
\qquad \Var{f} := \Vert f\Vert_2^2 - (
\overline{f})^2.
\]
The $L^2$ norm will be used in several different ways.
For a possible time-$t$ configuration $\mathbf{x}(t)$ of the averaging process,
the quantity $\Vert\mathbf{x}(t)\Vert_2$ is a number, and so the quantity
$\Vert\bX(t)\Vert_2$ appearing in the proposition below is a random variable.
The \emph{Dirichlet form} of the associated Markov chain is
\[
\EE(f,f) := {\frac{1}{2}} n^{-1} \sum
_i \sum_{j \ne i} (f_i
- f_j)^2 \nu_{ij} = n^{-1} \sum
_{\{i,j\}} (f_i - f_j)^2
\nu_{ij},
\]
where $\sum_{\{i,j\}}$ indicates summation over \emph{unordered} pairs.
The spectral gap of the chain, defined as the gap between eigenvalue $0$
and the second eigenvalue of $\NN$, is characterized as
%
%
\begin{equation}
\label{lambda-char} \lambda= \inf_f \biggl\{ \frac{\EE(f,f)}{\Var(f)} \dvt
\Var(f) \neq0 \biggr\}.
\end{equation}
The proof of the following result uses only the ingredients above and the
Averaging process dynamics.
\begin{Proposition}[(Global convergence theorem~\cite{aldous-lanoue})]
\label{PGcon}
From an initial configuration $\mathbf{x}(0) = (x_i)$ with average zero,
the time-$t$ configuration $\bX(t)$ of the averaging process satisfies
%
%
\begin{equation}
\label{Model-2-bd} \Ex\bigl\Vert\bX(t)\bigr\Vert_2 \leq\bigl\Vert\mathbf{x}(0)
\bigr\Vert_2 \exp( - \lambda t/4), \qquad 0 \le t < \infty,
\end{equation}
where $\lambda$ is the spectral gap of the associated MC.
\end{Proposition}

To introduce the next result,
thinking heuristically of the agents who agent $i$ most frequently
meets as the ``local'' agents for $i$, it is natural to guess that
the configuration of the averaging process might become
``locally smooth'' faster than the ``global smoothness'' rate
implied by Proposition~\ref{PGcon}.
In this context we may regard the Dirichlet form
$\EE(f,f)$
as measuring the
``local smoothness'', more accurately the local roughness, of a
function $f$, relative to the local structure of the particular meeting process.
The next result implicitly bounds $\Ex\EE(\bX(t), \bX(t)) $ at
finite times
by giving an explicit bound for the integral over $0 \le t < \infty$.
Note that, from the fact that the spectral gap is strictly positive, we
can see directly that $\Ex\EE(\bX(t), \bX(t)) \to0$ exponentially
fast as $t \to\infty$;
Proposition~\ref{P-smooth} is a complementary non-asymptotic result.
\begin{Proposition}
\label{P-smooth}
For the Averaging process with arbitrary initial configuration $\mathbf{x}(0)$,
\[
\Ex\int_0^\infty\EE \bigl(\bX(t), \bX(t) \bigr)
\,dt = 2 \Var{\mathbf{x}(0)}.
\]
\end{Proposition}
This looks slightly magical because the identity does not depend on the
particular rate matrix $\NN$, but of course the definition of
$\EE$ involves $\NN$.

The next result illustrates the following.
\begin{GP}
Notions of \emph{duality}
are one of the interesting and useful tools in classical IPS,
and equally so in the social dynamics models we are studying.
\end{GP}

There is a formal definition of duality in~\cite{Liggett}, Section~2.3,
but we find the concept easier to understand via examples.
The duality between the Voter model and coalescing chains (recalled in
Section~\ref{VMdual}) is the simplest and most striking example. The
relationship we give in Proposition~\ref{P:Adual} for the Averaging model
is perhaps more representative of the general style
of duality relationships.
We repeat the ``random penny''
argument used in Proposition~\ref{L1}.
Now there are two pennies, and at any meeting there are independent
decisions to hold or pass each penny.
Analyzing how the two pennies move leads to a certain coupling (\cite
{aldous-lanoue}, Sec.~2.4) of two (half-speed) associated Markov chains
$(Z_1(t),Z_2(t))$ -- the motion of the two pennies is coupled because,
at a meeting of their two owners, one or both may be passed.
The coupled process has the following relationship with the Averaging process.
\begin{Proposition}[(The duality relation for the Averaging process)]
\label{P:Adual}
For the averaging model $\bX(t)$ started from a configuration $\mathbf{x}(0)$
which is a probability distribution over agents,
and for each $t$,
\[
\Ex \bigl(X_i(t) X_j(t) \bigr) = \Pr
\bigl(Z_1(t) = i, Z_2(t) = j \bigr),
\]
where $(Z_1(t),Z_2(t))$ denotes the coupled process started from random
agents $(Z_1(0),Z_2(0))$ chosen independently from $\mathbf{x}(0)$.
\end{Proposition}

On a simple geometry where one can calculate
$\Ex X_i(t)$ via Proposition~\ref{L1},
one could seek to use Proposition~\ref{P:Adual} to investigate
second-order properties of the Averaging process, and some open
problems of this type are mentioned in~\cite{aldous-lanoue}.

Finally, parallel to Proposition~\ref{PGcon} which builds upon a standard
$L^2$ bound for reversible Markov chains is the following result
\cite{aldous-lanoue}, inspired by the theory surrounding
log-Sobolev inequalities
\cite{Diaconis}.
For a configuration $\mathbf{x}$ which is a probability distribution
write
\[
\Ent(\mathbf{x}) : = - \sum_i x_i
\log x_i
\]
for the entropy of the configuration.
Consider the averaging process where the initial configuration is a
probability distribution.
By concavity of the function $- x \log x$ it is clear that in the
averaging process $\Ent(\bX(t))$ can only increase, and hence
$\Ent(\bX(t)) \uparrow\log n$ a.s.
(recall $\log n$ is the entropy of the uniform distribution).
\begin{Proposition}
\label{Pentropy}
For the Averaging process whose initial configuration is a
probability distribution $\mathbf{x}(0)$,\vspace*{-1pt}
\[
\Ex \bigl(\log n - \Ent \bigl(\bX(t) \bigr) \bigr) \le \bigl(\log n - \Ent
\bigl( \mathbf{x}(0) \bigr) \bigr) \exp(-\alpha t/2),
\]
where $\alpha$ is the log-Sobolev constant of the associated Markov chain.
\end{Proposition}

\section{The Voter model}
\label{sec:VM}

This model has a random update rule, most elegantly
implemented via a
``directed'' convention in the meeting model:
when agents $i,j$ meet, choose a (uniform) random direction and
indicate it using an arrow $i \to j$ or $j \to i$.

\begin{quote}
\begin{Votermodel*}
Initially each agent has a different ``opinion'' -- agent $i$ has
opinion~$i$.
When $i$ and $j$ meet at time $t$ with direction $i \to j$,
then agent $j$ adopts the current opinion of agent $i$.
\end{Votermodel*}
\end{quote}

This model and many variants have acquired a substantial literature in
the finite setting (as well as the classical infinite lattice setting)
so what's written here is far from complete.
In the classical setting, one traditionally assumed only two different
initial opinions.
In our finite-agent case, it seems more natural to take the initial
opinions to be all different.

For notation, write
\[
\VV_i(t) := \mbox{the set of $j$ who have opinion $i$ at time $t$}.
\]
So $\{ \VV_i(t), i \in\Agents\}$ is a random partition of $\Agents$.
Note that $\VV_i(t)$ may be empty, or may be non-empty but not contain $i$.
The number of different remaining opinions can only decrease with time.
Ultimately the process must be absorbed in one of the $n$
``everyone has same opinion'' configurations.
This fits one of our Section~\ref{sec:SCC} general possibilities for a
FMIE process:
absorption in a random one (of a small number of) ``ordered'' configurations.
A natural quantity of interest is this \textbf{consensus time}
\[
T^{\mathrm{voter}}: = \min \bigl\{t: \VV_i(t) = \Agents\mbox{ for
some } i \bigr\}.
\]

\subsection{Duality with coalescing chains}
\label{VMdual}

\begin{quote}
\begin{Coalescingchainsmodel*}
Initially each agent has a token -- agent $i$ has token $i$.
At time $t$ each agent $i$ has a (maybe empty) collection (cluster)
$\CC_i(t)$ of tokens.
When $i$ and $j$ meet at time $t$ with direction $i \to j$,
then agent $i$ gives his tokens to agent $j$;
that is,
\[
\CC_j(t+) = \CC_j(t-) \cup\CC_i(t-),\qquad
\CC_i(t+) = \varnothing.
\]
\end{Coalescingchainsmodel*}
\end{quote}
So $\{ \CC_i(t), i \in\Agents\}$ is a random partition of $\Agents$.
A natural quantity of interest is the \textbf{coalescence time}
\[
T^{\mathrm{coal}}: = \min \bigl\{t: \CC_i(t) = \Agents\mbox{ for
some } i \bigr\}.
\]
Regarding each non-empty cluster as a particle, the particles move as the
associated Markov chain at half-speed
(rates $\nu_{ij}/2$), moving independently until two particles meet
and thereby coalesce.
Note this factor $1/2$ in this section.

\begin{Proposition}[(The duality relation for the Voter model)]
For fixed $t$, $\{\VV_i(t), i \in\Agents\} \ed  \{\CC_i(t), i \in
\Agents\} $.
In particular
$T^{\mathrm{voter}}\ed T^{\mathrm{coal}}$.
\end{Proposition}
It is important to realize that
they are different as processes.
For fixed $i$, note that $|\VV_i(t)|$ can only change by $\pm1$, but
$|\CC_i(t)|$ jumps to and from $0$.

The proof in the classical lattice setting~\cite{MR538077} (which
refers to coalescing \emph{random walks} extends unchanged to our
setting of general symmetric rates $(\nu_{ij})$. One simply considers
a realization of the (directed) meeting process on a fixed time
interval $[0,t]$, uses this to construct the
coalescing chains process at time $t$, then observes that reversing the
direction of time and the arrows gives a construction of the Voter model
at time $t$.

Almost all of the literature on the finite Voter model has focussed on
estimating
$T^{\mathrm{voter}}\ed T^{\mathrm{coal}}$ in terms of the geometry,
and I will show some of
this work in Sections~\ref{sec:VMCG}--\ref{sec:King2}.
But there are other questions one can ask about the finite-time behavior,
to be mentioned in Section~\ref{sec:pre-asy}.

\subsection{The Voter model on the complete graph}
\label{sec:VMCG}

There are two ways to analyze $T^{\mathrm{voter}}_n$ on the complete
graph, both
providing some bounds on other geometries.
The first uses coalescing chains.
One aspect of \textbf{Kingman's coalescent}~\cite{beres} is the
continuous-time MC on states
$\{1, 2, 3, \ldots\}$ with standardized rates $\lambda_{k,k-1} = {k
\choose2}, k \geq2$.
For that chain, the mean hitting time to state $1$ is
\[
\Ex_m T^{\mathrm{hit}}_1 = \sum
_{k=2}^m 1\Big/ \pmatrix{k
\cr
2} = 2 \biggl(1 - {
\frac{1}{m}} \biggr)
\]
and in particular $\lim_{m \to\infty} \Ex_m T^{\mathrm{hit}}_1 = 2$.
Now when we consider coalescing chains on the complete $n$-graph, the
number of clusters evolves as the
continuous-time MC on states
$\{1, 2, 3, \ldots, n\}$ with rates $\lambda_{k,k-1} = {\frac
{1}{n-1}} {k \choose2}$.
So
$\Ex T^{\mathrm{coal}}_n = (n-1) \times2 (1 - {\frac{1}{n}}) $
and in particular
%
%
\begin{equation}
\label{voter-Kn} \Ex T^{\mathrm{voter}}_n = \Ex T^{\mathrm{coal}}_n
\sim2 n.
\end{equation}

The second way is to consider the Voter model directly, specifically
the variant of the voter model with only 2 opinions. Consider
the number $X(t)$ of agents with the first opinion.
On the complete $n$-graph, $X(t)$ evolves as the continuous-time MC on states
$\{0,1,2,\ldots,n\}$ with rates
\[
\lambda_{k,k+1} = \lambda_{k,k-1} = {\frac{k(n-k)}{2(n-1)}}.
\]
This process arises in classical applied probability (e.g., as the
Moran model in population genetics).
We want to study the absorption time
\[
T^{\mathrm{hit}}_{0,n} := \min \bigl\{t: X(t) = 0 \mbox{ or } n \bigr
\}.
\]
By general birth-and-death formulas, or by comparison with simple RW,
\[
\Ex_k T^{\mathrm{hit}}_{0,n} = {\frac{2(n-1)}{n}} \bigl( k
(h_{n-1} - h_{k+1}) + (n-k) (h_{n-1} -
h_{n-k+1}) \bigr) \qquad\mbox{where } h_m := \sum
_{i=1}^m 1/i.
\]
This is maximized by $k = \lfloor n/2 \rfloor$, and
\[
\max_k \Ex_k T^{\mathrm{hit}}_{0,n} \sim(2
\log2) n.
\]

Now we can couple the original Voter model ($n$ different initial
opinions) with the variant with only 2 opinions, initially held by
$k$ and $n-k$ agents.
(Just randomly assign these two opinions, initially.)
From this coupling, we see
\begin{eqnarray*}
\Pr_k \bigl(T^{\mathrm{hit}}_{0,n} > t \bigr) &\le&\Pr
\bigl(T^{\mathrm{voter}}_n > t \bigr),
\\
\Pr_k \bigl(T^{\mathrm{hit}}_{0,n} > t \bigr) &\ge&{
\frac{2k(n-k-1)}{n(n-1)}} \Pr \bigl(T^{\mathrm{voter}}_n > t \bigr),
\end{eqnarray*}
where the fraction is the chance that two specified agents are assigned
different opinions, initially.
In particular, the latter with $k = \lfloor n/2 \rfloor$ implies
\[
\Ex T^{\mathrm{voter}}_n \le \bigl(4 \log2 + o(1) \bigr) n.
\]
Comparing with the correct asymptotics (\ref{voter-Kn}), we see that
although constant is not precise, the order of magnitude is correct.

\subsection{A crude bound for general geometries}

Consider a general geometry $(\nu_{ij})$, and write
$\nu(A,A^c) := \sum_{i \in A,j\in A^c} \nu_{ij}$.
Suppose the flow rates satisfy, for some constant $\kappa> 0$,
\[
\nu \bigl(A,A^c \bigr) := \sum_{i \in A, j \in A^c}
n^{-1} \nu_{ij} \ge \kappa\frac{|A| (n - |A|)}{n(n-1)}.
\]
The right side is normalized so that on the complete graph this holds
with $\kappa= 1$.
We can repeat the analysis above -- the process $X(t)$ now moves at least
$\kappa$ times as fast as on the complete graph, and so
\[
\Ex T^{\mathrm{voter}}_n \le \bigl(4 \log2 + o(1) \bigr) n/ \kappa.
\]
This style of argument can be used for many FMIE processes.
\begin{GP}
Bottleneck statistics give crude general
bounds.
\end{GP}

Often the speed at which the asymptotic behavior sets in depends on
some quantification of the connectivity of the geometry.
If for some $m$ the quantity
\[
\phi(m) = \min \bigl\{\nu \bigl(A,A^c \bigr) : |A| = m \bigr\},
\qquad1 \le m \le n-1
\]
is small, it indicates a possible ``bottleneck'' subset of size $m$.
It is convenient to combine the family $(\phi(m), 1 \le m \le n/2)$
into a single statistic, but the appropriate way to do this is
(FMIE) model-dependent.
In the Voter model case above, we used
\[
\kappa: = \min_A \frac{n (n-1) \nu(A,A^c)}{|A| (n - |A|)} = n(n-1) \min_m
\frac{\phi(m)}{m(n-m)}.
\]
Quantities like $\kappa$ are descriptive statistics of a weighted graph,
often given names like ``isoperimetric constants'', but
``bottleneck statistics'' seems a more informative name.

The best known use of bottleneck statistics is Cheeger's inequality
(which bounds the spectral gap of a Markov chain in terms of $\kappa$)
and its generalizations using other bottleneck statistics -- see
\cite{MR2341319} for an detailed treatment.
Such techniques are used extensively in the theory of algorithms, where
$\kappa$ is called
\emph{conductance} -- see, e.g.,~\cite{MR2809775} for analysis of a
gossip algorithm.

\subsection{Coalescing chains on general geometries}

Estimating $T^{\mathrm{coal}}$ is clearly related to study of the
\emph{meeting time}
$T^{\mathrm{meet}}_{ij}$ of two independent copies of the chain
started at $i$
and $j$, a topic that arises in other contexts.
Recall that on the complete graph the mean coalescence time is asymptotically
twice the mean meeting time. Also note that
under enough symmetry (e.g., continuous-time random walk on the
discrete torus)
the relative displacement between the two copies evolves as the same
random walk run at twice the speed,
and study of $T^{\mathrm{meet}}_{ij}$ reduces to study of $T^{\mathrm
{hit}}_k$.
All this suggests there is a close connection between the coalescence
time and hitting times for the associated chain.
The following result was conjectured long ago but only recently proved.
For any chain, define a statistic
\[
\tau^* := \max_{i,j} \Ex_i T^{\mathrm{hit}}_j
.
\]
\begin{Theorem}[(\cite{oliveira_coal})]
\label{T_oliveira}
There exist numerical constants $C_1, C_2 < \infty$ such that, for any
finite irreducible reversible Markov chain,
\[
\max_{i,j} \Ex T^{\mathrm{meet}}_{ij} \le C_1
\tau^* \quad\mbox{and}\quad \Ex T^{\mathrm{coal}}\le C_2 \tau^*.
\]
\end{Theorem}

The proof involves intricate analysis that I will not try to outline here.
If the stationary distribution $\pi$ contains exponentially rare states
then the theorem is uninformative. But in our setting of FMIE
processes, the associated chain always has \emph{uniform} stationary
distribution.
One can create examples of geometries in which the bound is not the correct
order of magnitude, but I suspect that in any ``natural'' geometry
Theorem~\ref{T_oliveira} does give the correct order.

A different style of general result given recently in
\cite{cooper_CRW} shows that on a regular graph
$\Ex T^{\mathrm{coal}}= O(n/\lambda)$
where $\lambda$ is the spectral gap of the associated MC.
But this bound is of correct order only for expander-like graphs.

\subsection{How general is the Kingman coalescent calculation?}
\label{sec:King2}

Return to the setting of Section~\ref{VMdual} -- chains are associated
chains run at half speed.
Write $\taumeet$ for the mean meeting time from independent uniform starts.
In a sequence of chains with $n \to\infty$, impose a condition such
as the following.
For each $\eps> 0$
%
%
\begin{equation}
\label{VMau} n^{-2} \bigl| \bigl\{(i,j): \Ex T^{\mathrm{meet}}_{ij}
\notin(1 \pm\eps) \taumeet \bigr\} \bigr| \to0 \qquad\mbox{as } n \to\infty.
\end{equation}
%

\begin{OP}
Assuming (\ref{VMau}), under what further conditions can one prove
%
%
\begin{equation}
\Ex T^{\mathrm{coal}}\sim2 \taumeet\qquad\mbox{as $n \to\infty$?}
\label{OPKingman}
\end{equation}
\end{OP}

That such a result is anticipated has been folklore since 1989, when
Cox~\cite{cox}
proved (\ref{OPKingman}) for the torus $[0,m-1]^d$ in dimension $d \ge2$.
(In that geometry, $\taumeet\sim m^d R_d$ for $d \ge3$,
where $R_d$ is the mean number of visits of standard random walk on
$\Ints^d$ to the origin.)
See~\cite{blythe} for a physicist's view.

One could try to prove this in two stages.

\textbf{Stage 1.}
For fixed $m$, show that the mean time for $m$ initially independent
uniform walkers to coalesce should be
$\sim2(1 - {\frac{1}{m}}) \taumeet$.

\textbf{Stage 2.}
Show that for $m(n) \to\infty$ slowly, the time for the initial $n$
walkers to
coalesce into $m(n)$ clusters is $o(\taumeet)$.

It is likely that Stage 1 can be proved without further assumptions,
as follows.
From old results on mixing times~\cite{aldous-oldmixing}, condition
(\ref{VMau})
is enough to show that the variation distance mixing time $\taumix=
o(\taumeet)$. So -- as a prototype use of
$\taumix$ -- by considering time intervals of length $\tau$,
for $\taumix\ll\tau\ll\taumeet$, the events
``a particular pair of walkers meets in the next $\tau$-interval''
are approximately independent.
This makes the ``number of clusters'' process behave (to first order)
as the Kingman coalescent as the number decreases from $m$ to $1$.
This approach is implemented in~\cite{MR2114986} for $\Ints^2_m$, but
it should hold for general geometries.

On the other hand, Stage 2 requires some different assumptions to
control short-time
behavior.
See~\cite{Oliveira12} for a recent extensive study of this problem, and
the current best results.

\subsection{Pre-asymptotic behavior}
\label{sec:pre-asy}

As mentioned earlier, almost all of the literature on the finite Voter
model has focussed on estimating
$T^{\mathrm{voter}}$ in terms of the geometry.
But there are other aspects of finite-time behavior one might study, and
we can uncover some such aspects by
comparison with what we have studied in the Averaging process.
Recall we start with all opinions different.

\textbf{1.} If the proportions of agents with the various opinions are
written as $\mathbf{x}= (x_i)$,
the statistic $q := \sum_i x_i^2$ is one measure of concentration
versus diversity of opinion.
This suggests studying $Q(t) := \sum_i ( n^{-1} |\VV_i(t)|)^2$.
Duality implies
\[
\Ex Q(t) = \Pr \bigl( T^{\mathrm{meet}}\le t \bigr),
\]
where $T^{\mathrm{meet}}$ is the meeting time for independent chains
with uniform starts.
One can study the distribution of $T^{\mathrm{meet}}$ in special geometries,
and more generally can often apply the explicit Exponential distribution
approximation from~\cite{aldous-fill}, Prop.~3.23.

\textbf{2.} A corresponding ``local'' measure of concentration versus diversity
is the probability that agents $(I,J)$ chosen with probability $\propto
\nu_{ij}$
(``random neighbors'') have same opinion at time~$t$.
Refinements of this ``clustering'' behavior have been studied in the
classical two-dimensional lattice setting~\cite{cox-griffeath}.

\textbf{3.} The diversity statistic $q := \sum_i x_i^2$ emphasizes
large clusters (large time); the
statistic $\Ent(\mathbf{x}) = - \sum_i x_i \log x_i$ emphasizes
small clusters (small time).
So one could consider
\[
\mathrm{E}(t) : = - \sum_i \bigl(
n^{-1} \bigl|\VV_i(t)\bigr| \bigr) \log \bigl( n^{-1}\bigl |
\VV_i(t)\bigr| \bigr).
\]
This has apparently not been studied.

\subsection{Approximating finite graphs by infinite graphs}
\label{sec:infinite}

\begin{GP}
Certain special families of geometries have $n \to\infty$ local weak limits
(which are infinite rooted random networks), and aspects of FMIE models
on the finite geometries can often be related to the behavior of the
model on the infinite geometry.\looseness=1
\end{GP}

Two familiar examples of such geometries are
\begin{longlist}[(ii)]
\item[(i)]For the torus $\Ints^d_m$, the $m \to\infty$ limit is the
infinite lattice $\Ints^d$;
\item[(ii)]For the ``random graphs with prescribed degree
distribution'' geometry,
the limit is a Galton-Watson tree with a size-biased offspring
distribution after the first generation;
\end{longlist}
and two less familiar examples are
\begin{longlist}[(iii)]
\item[(iii)]For the complete graph with random edge-weights, the limit
is the
PWIT~\cite{aldous-steele}, an infinite-degree tree where the weights
on the edges at a given vertex form a rate-$1$ Poisson process on
$[0,\infty)$;
\item[(iv)]A more elaborate graph~\cite{aldous-complex} built over (iii),
designed to have a
$n \to\infty$ local weak limit which is more complex (not just a
tree) yet still fairly tractable.
\end{longlist}

The relation between finite and infinite geometries is simplest for the
``epidemic'' (FPP) type models later, but
also can be used for Markov chain related models, starting from the following.

\begin{Localtransienceprinciple*}
For a large finite-state Markov chain whose behavior near a state $i$
can be approximated be a transient infinite-state chain, we have
\[
\Ex_\pi T^{\mathrm{hit}}_i \approx R_i/
\pi_i,
\]
where $R_i$ is defined in terms of the approximating infinite-state
chain as
$\int_0^\infty p_{ii}(t) \,dt$.
The approximation (an instance of the Poisson clumping heristic \cite
{aldous-PCH}) comes from the finite-state mean hitting time formula
(\cite{aldous-fill}, Lemma~2.11)
via a ``interchange of limits'' procedure which requires ad hoc justification.
\end{Localtransienceprinciple*}

In the context of the Voter model, if there is a local weak limit of
the finitie geometries,
the local transience principle implies that the property
$\Ex T^{\mathrm{voter}}= \Theta(n)$ is essentially equivalent to
transience of
the chain
on the limit infinite geometry.

\section{Two further models}

\subsection{The Compulsive Gambler process}

\begin{quote}
\begin{TheCompulsiveGamblerprocess*}
Initially each agent has some (nonnegative real-valued) amount of money.
Whenever two agents meet, they instantly play a fair game in which one
agent acquires the combined money. In other words,
if one has $a$ and the other has $b$ then the first acquires all $a+b$ with
chance $a/(a+b)$.
\end{TheCompulsiveGamblerprocess*}
\end{quote}
Note that on the complete graph geometry,
this process is just an augmentation of the
Kingman coalescent process.
On a general geometry, a configuration in which the set of agents with
non-zero money forms an ``independent set''
(no two are adjacent in the weighted graph) is obviously an absorbing
configuration, and conversely.

This process has apparently not been studied.
I mention it because it provides a simple example of the ``disordered'' limit
behavior mentioned in Section~\ref{sec:SCC}.

\subsection{The Interchange (Exclusion) process}

\begin{quote}
\begin{TheInterchangeprocess*}
Initially each agent has a different token.
When two agents meet, they exchange tokens.
\end{TheInterchangeprocess*}
\end{quote}
I will call this the Interchange process; the ``two type'' version
is the classical Exclusion process.
Though this model is not particularly natural in the social networks context,
the result below (which established a longstanding conjecture)
illustrates a result for which our general ``weighted graphs'' setting
is the
(mathematically) correct setting.

The motion of an individual token is as the associated Markov chain, so
has some
spectral gap $\lambda_{\mathrm{MC}}(\NN) > 0$.
The whole Interchange process is itself a reversible Markov chain, so
has a spectral gap $\lambda_{\mathrm{IP}}(\NN)$ which by a general
``contraction'' fact (\cite{aldous-fill}, Sec.~4.6.1) satisfies
$\lambda_{\mathrm{IP}}(\NN) \leq\lambda_{\mathrm{MC}}(\NN) $.
\begin{Proposition}[(\cite{caputo})]
For every geometry $\NN= (\nu_{ij})$ we have $\lambda_{\mathrm
{IP}}(\NN) = \lambda_{\mathrm{MC}}(\NN) $.
\end{Proposition}
The proof involves intricate analysis that I will not try to outline here.

\section{The Pandemic process}
\label{sec:Pandemic}

Recall the discussion in Section~\ref{sec:Pan}.
For the purposes of this section, the Pandemic process is the same as
the FPP
process; we will use the language of ``infectives'' and ``time''
(rather than ``percolation distance'') and write $T^{\mathrm
{fpp}}_{ij}$ for the
time taken for an epidemic started at $i$ to reach $j$.

In this section, we will recall the well-understood behaviors on the
lattice (Section~\ref{sec:FPP_torus}) and the complete graph (Section~\ref{sec:FPP_complete}).
Each type of behavior raises questions about the extent to which
similar behavior remains true in more general geometries
(Sections~\ref{sec:WLLN} and~\ref{sec:FPP_other}) which have not been
systematically studied.

\subsection{FPP on the torus}
\label{sec:FPP_torus}

We quote the classical result for FPP on the infinite lattice,
with brief remarks.
\begin{Theorem}[(Shape theorem for lattice FPP~\cite{kesten})]
Consider the FPP process on the edges of $\Ints^d$
with $\operatorname{Exponential}(1)$ edge-lengths, starting from the origin.
Write $\SS(t) \subset\Ints^d$ for the infected set at time $t$,
and $\bar{\SS}(t) \subset\Reals^d$ for its fattening (replace each
point by a cube of side $1$).
There is a non-random closed convex set $B = B_d$ such that,
for each $0 < \eps< 1$,
\[
\Pr \bigl( (1 - \eps) t B \subseteq\bar{\SS}(t) \subseteq(1+\eps) t B \bigr)
\to1 .
\]
\end{Theorem}
\begin{itemize}
\item The key proof ingredient is the \emph{subadditive ergodic theorem}.
\item The result hold for more general IID edge-times.
\item$B_d$ is \emph{not} the unit ball in Euclidean norm.
\item For later use write $b_d$ for the ($d$-dimensional) ``volume'' of $B_d$.
\item Understanding the asymptotic variance of $T^{\mathrm{fpp}}_{ij}$
is a
famous hard problem; see, e.g.,~\cite{benaim}.
\end{itemize}
Consider now the FPP process on the discrete torus $\Ints^d_m$.
We use our standardization convention (\ref{standardization}),
so that rates $\nu_{ij} = 1/(2d)$ for each edge $(i,j)$.
The following three ``first-order'' results are immediate from the
shape theorem.
Asymptotics are as $m \to\infty$ for fixed $d$.
%
%
\begin{equation}
\Ex\bigl|\SS(t)\bigr| \sim(2d)^{-1} b_d t^d \mbox{ over
} 1 \ll t \ll m.
\end{equation}
%
The same holds (in $L^1$) for $|\SS(t)|$ itself.
Second,
the time $T^{\mathrm{pan}}_*(m)$ until all agents are infected satisfies
%
%
\begin{equation}
m^{-1} T^{\mathrm{pan}}_*(m) \to c_d \qquad\mbox{a.s.}
\end{equation}
for a constant $c_d$ defined in terms of the shape of $B_d$.
Third, the proportion $X_m(t)$ of infected agents satisfies
%
%
\begin{equation}
\sup_{0 \le s < \infty} \bigl|X_m(ms) - F_d(s)\bigr|
\to_p 0
\end{equation}
for a certain function $F_d$ defined in terms of the shape of $B_d$.

\subsection{Toward a general WLLN}
\label{sec:WLLN}

Having a limit ``shape'' represented by $B_d$ is very special
to the lattice setting.
What about the general setting of rates $(\nu_{ij})$?
Suppose we have a result of the kind
\[
\mbox{(*)\quad provided $i$ and $j$ are not close, then $T^{\mathrm{fpp}}_{ij}$
is close to its mean.}
\]
Then the random set $\SS_i(t)$ of infectives at $t$ from initial $i$
is approximately the deterministic set
$\{j: \Ex T^{\mathrm{fpp}}_{ij} \le t \}$.
So, if the essence of a ``shape theorem'' is a deterministic
approximation to the infected set, it is essentially just a WLLN
(weak law of large numbers).

A precise conjecture is formulated below.
Everything depends on $n$ (not written) and limits are as $n \to\infty$.
The conclusion we want is
%
%
\begin{equation}
\label{conj-1} \frac{T^{\mathrm{fpp}}_{ij}}{\Ex T^{\mathrm{fpp}}_{ij}} \to_p 1.
\end{equation}
We work in the FPP model, writing $\xi_{ab}$ for the time to traverse
edge $(a,b)$.

The obvious obstacle to the desired conclusion is that there may be a
set $A$
(with $i \in A, j \notin A$) such that, in the random percolation
path $\pi_{ij}$, the time taken to traverse the edge $A \to A^c$ is
not $o( \Ex T^{\mathrm{fpp}}_{ij})$.
One could conjecture this is the only obstacle; here is a slightly weaker
conjecture.
\begin{Conjecture}
\label{Conj-1}
With arbitrary rates $(\nu_{ij})$, if
%
%
\begin{equation}
\label{Tconc} \frac{\max\{\xi_{ab}: (a,b) \mbox{ edge in } \pi_{ij} \}} {
\Ex T^{\mathrm{fpp}}_{ij}} \to_p 0
\end{equation}
then (\ref{conj-1}) holds.
\end{Conjecture}

It is intuitively clear (and not hard to prove) that condition (\ref
{Tconc}) is necessary.
Proving the conjecture would hardly be a definitive result (because
condition (\ref{Tconc}) is not a condition directly on $(\nu_{ij}))$
but would
represent a start on understanding when the WLLN holds.

\subsection{FPP on the complete $n$-vertex graph}
\label{sec:FPP_complete}

Here we have rates
\[
\nu_{ij} = 1/(n-1),\qquad j \neq i
\]
and we start with one infective.
This is a fundamental stochastic process whose basic properties have
been rediscovered many times from different viewpoints, but there seems
no authoritative survey.
Our focus is on the randomly-shifted logistic limit result,
Proposition~\ref{Plogistic}. We give an outline
taken from~\cite{aldous-gossip} (a similar treatment is in \cite
{MR2439767}), though essentially the same
result can be found in older genetics literature.
There are two reasons for this particular focus. From this one result,
we can read off many results like (\ref{DGG},~\ref{DGL}), a fact
which has not been appreciated in the combinatorics-inspired literature.
And secondly, it plays a central role is investigating FMIE processes
build over Pandemic, as we shall see in
Section~\ref{sec:built_over_P}.

Write $D(k) = D_n(k) = $ time at which $k$ people are infected.
The observation
%
%
\begin{equation}
\label{knk} \mbox{the r.v.'s $D(k+1)- D(k)$ are independent with $
\operatorname{Exponential}\bigl({\frac{k(n-k)}{n-1}}\bigr)$ distribution}
\end{equation}
is a basic starting place for analysis.
For instance, the time $T^{\mathrm{pan}}_*$ until the whole population is
infected is just $D(n)$, and so
\[
\Ex T^{\mathrm{pan}}_* = \sum_{k=1}^{n-1}
\Ex \bigl(D(k+1) - D(k) \bigr) = (n-1) \sum_{k=1}^{n-1}{
\frac{1}{k(n-k)}} \sim2 \log n.
\]

And the time $T^{\mathrm{fpp}}_{{\mathrm{rand}}}$ until a random
person (amongst the $n-1$ initially uninfected) is infected is
\begin{eqnarray}
\label{Euler} \Ex T^{\mathrm{fpp}}_{{\mathrm{rand}}} &=& \Ex D(U),\qquad\mbox{$U$
uniform on } \{2,3,\ldots,n\}
\nonumber
\\
&=& {\frac{1}{n-1}} \sum_{k=2}^n \Ex
D(k)
\nonumber
= {\frac{1}{n-1}} \sum_{k=1}^{n-1}
(n-k) \Ex \bigl(D(k+1) - D(k) \bigr)
\nonumber
\\
&=& \sum_{k=1}^{n-1} {\frac{1}{k}}
\sim\log n.
\end{eqnarray}
This analysis can be pursued to get distributional results, but let us first
look at different methods that indicate the behavior without extensive
calculation.
The idea is to analyze separately the \emph{initial phase} when $o(n)$
people are infected,
and the \emph{pandemic phase} when $\Theta(n)$ people are infected.

The ``deterministic, continuous'' analog of our ``stochastic,
discrete'' model of an epidemic is the
\emph{logistic equation}
\[
F^\prime(t) = F(t) \bigl(1 - F(t) \bigr)
\]
for the proportion $F(t)$ of a population infected at time $t$.
A solution is a shift of the basic solution
\[
F(t) = {\frac{e^t}{1+e^t}},\qquad - \infty< t < \infty \qquad(\mbox{the logistic
function}).
\]
Setting $M_n(t)$ for the number of individuals infected in our
stochastic process,
then scaling to the proportion
$X_n(t) = M_n(t)/n$, we find
\[
\Ex \bigl(dX_n(t)|\FF_n(t) \bigr) = {\frac{n}{n-1}}
X_n(t) \bigl(1 - X_n(t) \bigr) \,dt.
\]
By including variance estimates, one could prove the following
formalization of the idea that, during
the pandemic phase, $X_n(t)$ behaves as $F(t)$ to first order.
Precisely,
for fixed $\eps$ and $t_0$ we have (in probability)
%
%
\begin{equation}
\label{XDFF} \bigl(X_n \bigl(D_n(n \eps) +t \bigr), 0
\le t \le t_0 \bigr) \to \bigl(F \bigl( F^{-1}(\eps) + t
\bigr), 0 \le t \le t_0 \bigr).
\end{equation}
That is, we ``start the clock'' when a proportion $\eps$ of the
population is infected.
This can now be reformulated more cleanly in terms of the time $G_n =
D_n(n/2)$ at which half the population is infected:
\[
\sup_{-G_n \le t < \infty} | n^{-1}M_n(G_n + t) -
F(t)| \to0 \qquad\mbox{in probability}.
\]
For later use, consider the times at which each process equals
$k_n/n$, where $k_n/n \to0 $ slowly. We find
\[
G_n + \log(k_n/n) = D_n(k_n) +
o(1)
\]
which rearranges to
%
%
\begin{equation}
\label{Dn7} D_n(k_n) - \log k_n =
G_n - \log n + o(1).
\end{equation}

Now turn to the initial phase of the epidemic.
On a fixed initial time interval $[0,t_0]$, the process $M_n(t)$ of
number of infectives converges in distribution
to the process $M_\infty(t)$ for which the times $D(k)$ satisfy
\[
\mbox{the r.v.'s $D(k+1) - D(k)$ are independent with $
\operatorname{Exponential}(k)$ distribution}
\]
and this is the classic \emph{Yule process}, for which it is well
known that
%
%
\begin{eqnarray}
\label{Yule_lim} &M_\infty(t) \mbox{ has Geometric $
\bigl(e^{-t} \bigr)$ distribution}&
\nonumber
\\
&e^{-t} M_\infty(t) \to\EE_\infty\mbox{ a.s. as } t
\to\infty, \mbox{ where } \EE_\infty\mbox{ has $\operatorname{Exponential}(1)$
distribution.}&
\end{eqnarray}
Now calculate informally, for $1 \ll k \ll n$,
\begin{eqnarray*}
\Pr \bigl(D_n(k) \le t \bigr) &=& \Pr \bigl(M_n(t) \ge k
\bigr) \approx\Pr \bigl(M_\infty(t) \ge k \bigr)
\\
&=& \Pr \bigl(e^{-t} M_\infty(t) \ge k e^{-t} \bigr)
\approx\Pr \bigl(\EE_\infty\ge ke^{-t} \bigr) = \exp
\bigl(-ke^{-t} \bigr).
\end{eqnarray*}
In other words
%
%
\begin{equation}
\label{Dnk} D_n(k) - \log k \approx_d G, \qquad1 \ll k
\ll n,
\end{equation}
where $G$ has the Gumbel distribution
$\Pr(G \le x) = \exp( - e^{-x}) $.
Comparing with (\ref{Dn7}) gives (\ref{GnG}) below.
\begin{Proposition}[(The randomly-shifted logistic limit)]
\label{Plogistic}
For the simple epidemic process on the complete $n$-vertex graph,
there exist random $G_n$ such that
%
%
\begin{eqnarray}
\label{GnG} \sup_{-G_n \le t < \infty} \bigl| n^{-1}M_n(G_n
+ t) - F(t)\bigr| &\to&0 \qquad\mbox{in probability},
\nonumber
\\
G_n - \log n &\cd& G,
\end{eqnarray}
where $F$ is the logistic function and $G$ has Gumbel distribution.
\end{Proposition}

One could prove this by formalizing the arguments above, but there is a
more efficient
though less illuminating way.
The result is essentially equivalent to the assertion
%
%
\begin{equation}
\label{Dnun} D_n \bigl(\lfloor un\rfloor \bigr) - \log n \cd
F^{-1}(u) + G, \qquad 0<u<1.
\end{equation}
The basic observation (\ref{knk}) allows one to write down an explicit
expression for the Fourier transform of
the left side, and one just needs to work through the calculus to check
it converges to the Fourier transform of
the right side.
(Transform methods are of course very classical in applied probability,
and their use in epidemic models goes back at least to~\cite{bailey}.)

The same transform argument (then taking a different limit regime) can
be used to prove the following result
(e.g.,~\cite{MR1965974}) for the time $D_n(n)$
until the entire population is infected:
%
%
\begin{equation}
\label{DGG} D_n(n) - 2 \log n \cd G_1 + G_2,
\end{equation}
where $G_1, G_2$ are independent with Gumbel distribution. While not
obvious from the statement of the
proposition above, it can in fact be deduced as follows. From (\ref{Dnun}),
the time $D_n(n/2)$ until half the population is infected satisfies
\[
D_n(n/2) - \log n \cd G.
\]
Then symmetry
(under $k \to n-k$) of the transition rates (\ref{knk}) implies that
$D_n(n)$ is essentially distributed as the sum of two independent
copies of
$D_n(n/2)$.

Reconsider the time $D_n(U)$ until a random person is infected.
Proposition~\ref{Plogistic} implies it has asymptotic distribution
%
%
\begin{equation}
\label{DGL} D_n(U) - \log n \cd G + L,
\end{equation}
where $L$ has logistic distribution function, independent of the Gumbel
time-shift $G$.
So we expect
\[
\Ex D_n(U) - \log n \to\Ex G + \Ex L.
\]
Now $\Ex L = 0$ by symmetry ($L \ed-L$) and one can calculate (or look up)
that $\Ex G = $ Euler's constant, so the limit theorem is consistent
with the exact formula
(\ref{Euler}).

\subsection{In what other geometries does similar behavior hold?}
\label{sec:FPP_other}

The extent to which different aspects of the ``mean-field'' behavior in
Section~\ref{sec:FPP_complete}
extend to more general geometries is not well understood.
For the ``random graphs with prescribed degree distributions'' (with
appropriate moments finite) geometry, one expects
to see qualitatively similar behavior (the pandemic phase to have duration
$O(1)$, centered at time $O(\log n)$), and almost precise analogs
have indeed been proved -- e.g.,~\cite{bhamidi1}, Theorem~3.1(b) gives
the analog
of (\ref{DGL}). There has been extensive subsequent work on related
problems and related ``locally tree-like''\vadjust{\goodbreak} geometries, but these proofs
are substantially
more technical than for the complete graph.
And even on the Hamming cube $\{0,1\}^d$ the behavior of FPP is not well
understood~\cite{MR1221168}.

To think about possible ``general theory'', let me suggest a loose
analogy with the
established theory of random walk on infinite graphs~\cite{woess}.
In that theory, most geometries fall into one of two categories,
exemplified by $\Ints^d$ and $d$-trees,
according to whether
\[
\Pr_i \bigl(X(t) = i \bigr) \to0 \qquad\mbox{as $ t^{-\alpha} $ or
$e^{-\beta t}$}.
\]
Properties involving spectral radius, isoperimetric inequalities, amenability,
(sub)linear rate of escape can then be related to this categorization.
That theory formalizes a distinction (informally, between ``finite
dimensional'' and ``infinite dimensional'') from the random walk viewpoint.

So by analogy, we have seen that the finite Pandemic process behaves
differently on $\Ints^d_m$ and on the complete graph.
Can we select one specific aspect of the different behavior as a useful
basis for classifying geometries?
Write $t_{ij} ;= \Ex T^{\mathrm{fpp}}_{ij}$.
One obvious distinction between these two geometries is that on $\Ints^d_m$ the means $t_{ij}$ vary (depending on distance between $i$ and
$j$) but on the complete graph, the $t_{ij}$ are the same
(to first order) for most pairs $(i,j)$.

So we can define a property (for a sequence of weighted graphs $(\nu_{ij})$)
which combines the idea of the WLLN (\ref{conj-1}) with the idea of the
$t_{ij}$ not varying much:
%
\begin{eqnarray}
\label{conj-infdim}&\mbox{there exist $t^{(n)}$ such that as $n \to
\infty$ for fixed $\eps> 0$}&
\nonumber
\\
&n^{-2} \bigl| \bigl\{(i,j): \Pr \bigl(T^{\mathrm{fpp}}_{ij}
\notin(1 \pm\eps ) t^{(n)} \bigr) > \eps \bigr\}\bigr | \to0. &
\end{eqnarray}
This property (which holds for the ``random graphs with prescribed
degree distribution'' geometry~\cite{bhamidi1}) is an attempt to
formalize the notion of a geometry being
``locally infinite-dimensional'' from the FPP viewpoint.
\begin{OP}
Study equivalences between this property and the
``geometry'' of $(\nu_{ij})$.
\end{OP}

Two brief comments.
\begin{longlist}[(ii)]
\item[(i)]We can't hope to establish (\ref{conj-infdim}) using only
isoperimetric inequalities -- they can only give order of
magnitude.
\item[(ii)]There is a loose analogy with the
cutoff phenomenon in finite Markov chains~\cite{persi-cutoff}, in that we
expect the property to hold in ``highly symmetric'' geometries.
But that analogy is discouraging in that there is no useful general
condition for cut-off.
\end{longlist}

\subsection{An epidemic process on the discrete torus $\Ints^2_m$
with short- and long-range interactions}

Consider the geometry $\NN$ with agents $\Ints^2_m$
and meeting rates
\begin{longlist}[(ii)]
\item[(i)]at rate $1$ an agent meets a uniform random neighbor;
\item[(ii)]at rate $m^{-\alpha}$ an agent meets a uniform non-neighbor.\vadjust{\goodbreak}
\end{longlist}
Here $0 < \alpha< 3$. Over this range the geometry qualitatively
interpolates between the complete-graph case and the nearest-neighbor
lattice case.
So we expect the center of the pandemic window to interpolate between orders
$\log m$ and $m$,
and the width of the pandemic window to interpolate between
orders $1$ and $m$.

For the Pandemic model over this geometry, it turns out that an analog
of the ``randomly-shifted logistic limit'' holds, albeit with a
different deterministic
function $F$
in place of the logistic, and different scaling. This is rather
counter-intuitive -- we envisage infectives as growing patches of mold in a Petri dish,
one might expect the detailed geometry of overlaps to play a role.
But the
key qualitative property is that, during the pandemic phase, the new
infectives arise from ``colonies'' started earlier in the pandemic
phase, not in the initial phase.
This enables us to derive the following equation
for the standardized function $F$ describing proportion of infectives
r over the pandemic window:
\[
1 - F(t) = \exp \biggl( - \int_{- \infty}^t F(s)
(t-s)^2 \,ds \biggr).
\]
See~\cite{aldous-gossip} for heuristics and~\cite{chat-gossip} for
careful treatment of a
closely related model.

\section{Some social analogs of epidemics}
\label{sec:built_over_P}

There is a large literature modeling spread of opinions, adoption of
innovations etc via social networks, and some of the models
fit our FMIE setting.
Instead of describing some standard model let me describe two new FMIE
models suggested by the ``social networks'' context.
Mathematically they are variants of Pandemic, and we will see
in Sections~\ref{sec:Deference}--\ref{sec:Fashionista} that on the
complete graph their behavior
can be derived from the
precise understanding of Pandemic given by Proposition~\ref{Plogistic}.

\subsection{The Deference model}
\label{sec:Deference}

\begin{quote}
\begin{TheDeferenceprocess*}
The agents are labelled $1$ through $n$, and
agent $i$ initially has opinion $i$.
When two agents meet, they adopt the same opinion, the smaller of the
two opinion-labels.
\end{TheDeferenceprocess*}
\end{quote}
Clearly opinion 1 spreads as Pandemic, so the ``ultimate'' behavior
of Deference is not a new question.

We will consider the complete graph geometry, and
study $(X^{n}_1(t), \ldots, X^{n}_k(t))$,
where $ X^{n}_k(t)$ is the proportion of the population with opinion
$k$ at time $t$.
We can rewrite Proposition~\ref{Plogistic} as
%
%
\begin{equation}
\label{XkF} \bigl(X^{n}_1(\log n + s), - \infty< s <
\infty \bigr) \cd \bigl(F(C_1 + s), - \infty< s < \infty \bigr),
\end{equation}
where $F$ is the logistic function and $C_1 = \log(\xi_1)$,
where $\xi_1 \ed \operatorname{Exponential}(1)$ arises as
the limit r.v. (\ref{Yule_lim}) associated with the Yule process.

The key insight is that opinions $1$ and $2$ and \ldots and $k$ combined
behave as one infection in Pandemic, except for a different initial phase,
and so the sum $\sum_{i=1}^k X^{(n)}_i(\cdot)$ behaves as in (\ref{XkF})
but with a different random shift $C_k$.
Because the initial phase is just a collection of $k$ independent Yule
processes,
we see that the sum has a representation of format (\ref{XkF}) with
%
%
\begin{equation}
\label{Cxi} C_k = \log( \xi_1 + \cdots+
\xi_k),\qquad k \ge1.
\end{equation}
Here the $(\xi_i)$ are i.i.d. $\operatorname{Exponential}(1)$.

So the $n \to\infty$ limit behavior of the Deference model in the
complete graph geometry is, over the pandemic phase,
%
%
\begin{eqnarray}
\label{XXk} && \bigl( \bigl(X^{n}_1(\log n + s),
X^{n}_2(\log n + s), \ldots, X^{n}_k(
\log n + s) \bigr), - \infty< s < \infty \bigr)
\\
&&\quad\to \bigl( \bigl(F(C_1 + s), F(C_2 + s) -
F(C_1 + s), \ldots, F(C_k + s) - F(C_{k-1} + s)
\bigr), - \infty< s < \infty \bigr).
\nonumber
\end{eqnarray}
%

\subsection{The Fashionista model}
\label{sec:Fashionista}

The Deference model envisages agents as
``slaves to authority''.
Here is a conceptually opposite ``slaves to fashion'' model, whose analysis
is mathematically surprisingly similar.

\begin{quote}
\begin{TheFashionistamodel*}
At the times of a rate-$\lambda$ Poisson process, a new fashion
originates with a uniform random agent, and is time-stamped.
When two agents meet, they each adopt the latest
(most recent time-stamp) fashion.
\end{TheFashionistamodel*}
\end{quote}
It is easy to see that, for the
random partition
of agents into ``same fashion'', there is a stationary distribution.

For the complete graph geometry, we can copy the analysis
of the Deference model.
Combining all the fashions appearing after a given time, these
behave as one infection in Pandemic
(over the pandemic window; newly-originating fashions have negligible effect),
hence their proportion behaves as a random time-shift of the logistic
curve $F$.
So when we study the vector
$(X^n_k(t), - \infty< k < \infty)$ of proportions of agents adopting
different fashions $k$,
we expect $n \to\infty$ limit behavior
of the form
%
%
\begin{equation}
\label{XXk2} \bigl(X^{n}_k(\log n + s), - \infty< k <
\infty \bigr) \cd \bigl( F(C_k + s) - F(C_{k-1} + s), -
\infty< k < \infty \bigr),
\end{equation}
where
$(C_k, - \infty< k<\infty)$
are the points
of some stationary process on $(- \infty, \infty)$.
Knowing this form for the $n \to\infty$ asymptotics, we can again
determine the distribution of $(C_i)$ by considering the initial
stage of spread of a new fashion. It turns out that
%
%
\begin{equation}
\label{Fash-limit} C_i = \log \biggl( \sum
_{j \le i} \exp( \gamma_j) \biggr) =
\gamma_i + \log \biggl( \sum_{k \ge1} \exp(
\gamma_{i - k} - \gamma_i) \biggr),
\end{equation}
where $(\gamma_j, - \infty< j < \infty)$ are the times of a
rate-$\lambda$ Poisson process on $(-\infty, \infty)$.
The second expression makes it clear that
$(C_i)$ is a stationary process.

Here is the outline argument for (\ref{Fash-limit}).
Consider the recent fashions at time $t = 0$ adopted by
small but non-negligible proportions\vadjust{\goodbreak} of the population.
More precisely, consider fashions originating during the time interval
$[- \log n + t_n, -\log n + 2t_n]$,
where $t_n \to\infty$ slowly.
For a fashion originating at time $- \log n + \eta$, the time-$0$ set of
adopting agents
will be a subset of the corresponding epidemic process, which we know
has proportional size $\xi\exp(- \eta) = \exp( - \eta+ \log\xi)$
where $\xi$ has $\operatorname{Exponential}(1)$
distribution.

The times $- \log n +\eta_j$ of origination of different fashions form
by assumption a rate-$\lambda$ Poisson process, and after we impose
i.i.d. shifts
$\log\xi_j$ we note (as an elementary property of Poisson processes)
that the shifted points $- \log n +\eta_j + \log\xi_j$ still form a
rate-$\lambda$ Poisson process, say
$\gamma_j$, on $(-\infty, \infty)$.
So the sizes of small recent fashion groups (that is letting $j \to-
\infty$),
for which overlap between fashions becomes negligible, are approximately
$\exp(\gamma_j)$.
Summing over $j \le i$ gives
\[
\sum_{j \le i} \exp( \gamma_j) \approx
F(C_i) \approx\exp(C_i)
\]
and we end up with the representation (\ref{Fash-limit}).

\subsection{Other geometries}

On the torus $\Ints^2_m$ the stationary distribution of Fashionista
is qualitatively similar to the \emph{dead leaves model}~\cite{bordenave}
and, as outlined in~\cite{aldous-FMIE}, one can apply simple
scaling arguments to analyze the diversity statistic
\[
s(m,\lambda) := \Pr(\mbox{two random agents share same fashion})
\]
and show
\[
s(m,\lambda) \sim c \lambda^{-2/3} m^{-2/3}\qquad \mbox{for }
m^{-1} \ll\lambda\ll m^2.
\]
The corresponding result on the complete graph is
\[
s(n,\lambda) \to s(\infty,\lambda) := \Ex\sum_k
\bigl( F \bigl(\lambda^{-1} C_k \bigr) - F \bigl(
\lambda^{-1}C_{k-1} \bigr) \bigr)^2
\]
and (\ref{Fash-limit}) provides a (complicated) explicit expression
for the right side.
In particular, we can see the $\lambda\to\infty$ behavior is
\[
s(\infty,\lambda) \sim c \lambda^{-1}
\]
for constant $c$.

Studying $(n,\lambda)$ in general geometries is an
\textbf{Open Topic}, which opens up other
questions to study about short-term behavior of
Pandemic.

\section{Other epidemic models}

There is a huge literature attempting more realistic epidemic
models, which I will not try to outline here.
The FMIE-type model that has been most extensively studied in the
mathematical literature is the Contact process (the SIS model
with Exponential infection times),
studied classically on $\Ints^d$~\cite{durrett-contact}.
A nice survey, in the style of this article, of recent rigorous work on
the Contact process on preferential attachment, small worlds and
power-law degree geometries can be found in
\cite{durrett-PNAS}, and the latest technical work in~\cite{MMVY12}.

\section*{Acknowledgements}

The current form of Conjecture~\ref{Conj-1}
emerged from conversations with Jian Ding.
And I thank an anonymous referee for many helpful expository suggestions.
Research supported by NSF Grant DMS-07-04159.


%

\end{document}